\documentclass[11pt,british,reqno]{amsart}

\usepackage[T1]{fontenc}

\usepackage{babel,eucal,url,amssymb,enumerate,amscd,color}

\theoremstyle{plain}
\newtheorem{lemma}{Lemma}[section]

\newtheorem{proposition}[lemma]{Proposition}

\newtheorem{theorem}[lemma]{Theorem}
\newtheorem{remark}[lemma]{Remark}

\newcommand{\Lie}[1]{\operatorname{\textsl{#1}}}

\newcommand{\Gtwo}{\ifmmode{{\rm G}_2}\else{${\rm G}_2$}\fi}

 \newcommand{\cyclic}{\mathop{\kern0.9ex{{+}\kern-2.2ex\raise-.28ex\hbox{\Large\hbox
 {$\circlearrowright$}}}}}

\def\sideremark#1{\ifvmode\leavevmode\fi\vadjust{\vbox to0pt{\vss
 \hbox to 0pt{\hskip\hsize\hskip1em
 \vbox{\hsize2.5cm\tiny\raggedright\pretolerance10000
 \noindent #1\hfill}\hss}\vbox to8pt{\vfil}\vss}}}%

\input xy
\xyoption{all}

\newfont{\eusm}{eusm10 scaled \magstep1}
\newfont{\eusmiii}{eusm10 scaled \magstep3}

\newcommand{\comp}{\makebox[7pt]{\raisebox{1.5pt}{\tiny $\circ$}}}

\title{Nearly K\"{a}hler homogeneous manifolds with positive curvature}
\author{J.~C.~Gonz{\'a}lez-D{\'a}vila}
\address[J.~C.~Gonz{\'a}lez~D{\'a}vila]{Department of Fundamental Mathematics\\
  University of La Laguna\\ 38200 La Laguna, Tenerife, Spain}
\email{jcgonza@ull.es}

\author{F.~Mart\'\i n~Cabrera}
\address[F.~Mart\'\i n~Cabrera]{Department of Fundamental Mathematics\\
  University of La Laguna\\ 38200 La Laguna, Tenerife, Spain}
\email{fmartin@ull.es}

\date{}

\date{\today}

\begin{document}

\maketitle

\begin{abstract}{\indent} We prove that a $2n$-dimensional compact homogeneous nearly K\"{a}hler
manifold with strictly positive sectional curvature is isometric to
${\mathbb C}P^{n},$ equipped with the symmetric Fubini-Study metric
or with the standard $\Lie{Sp}(m)$-homogeneous metric, $n =2m-1,$
or to $S^{6}$ as Riemannian manifold with constant sectional
curvature. This is a positive answer for a revised version of a conjecture given by Gray.

\vspace{4mm}

\vspace{2mm}

\end{abstract}





\section{Introduction}\indent
Berger \cite{Be} proved that a simply connected, even-dimensional,
normal homogeneous space of strictly positive sectional curvature
is  homeomorphic, in fact diffeomorphic, to a compact rank one
symmetric space. Here   we
firstly prove the following result:
\begin{theorem}\label{mean1} A simply connected,
$2n$-dimensional, normal homogeneous space of strictly positive
sectional curvature is isometric to a compact rank one symmetric
space or to the complex projective space
${\mathbb C}P^{n} =
\Lie{Sp}(m)/(\Lie{Sp}(m-1)\times \Lie{U}(1)),$ $n = 2m-1,$ $m\geq 2,$
equipped with the standard
$\Lie{Sp}(m)$-homogeneous Riemannian metric.
\end{theorem}
\noindent This is deduced as a
consequence of the Wallach's classification \cite{Wa} displayed in Section \ref{3symmetrcispaces} below.

Because the pinching constant $\delta =
\tfrac{\mathrm{min}\;K}{\mathrm{max}\;K}$ of the extremal values
of the sectional curvature for the standard
$\Lie{Sp}(m)$-homogeneous Riemannian metric is $\delta = \frac{1}{16}$ (see Proposition
\ref{ccp}), it follows that the metric on ${\mathbb C}P^{n}$
is not the symmetric Fubini-Study one.

On the other hand, Gray in \cite{G1} proposed the following
conjecture: {\em Let $(M,g,J)$ be a compact nearly K\"{a}hler manifold with strictly positive sectional
curvature. If the scalar curvature of
$(M,g)$ is constant, then $(M,g)$ is isometric to a complex projective space with a K\"ahler metric or to
 a $6$-dimensional sphere with a Riemannian
metric of constant sectional curvature.}

This conjecture is positive for K\"{a}hler manifolds  (see \cite{G3}). However, Sekigawa and Sato
  \cite{Se-Sa} showed that $\Lie{Sp}(2)/(\Lie{Sp}(1)\times \Lie{U}(1)),$ as
  a Riemannian $3$-symmetric space, is a counter example. Here we  prove

\begin{theorem}\label{mean2} If $M$ is a $2n$-dimensional compact homogeneous nearly K\"{a}hler manifold with strictly positive sectional
curvature, then $M$ is isometric to
\begin{enumerate}
\item[{\rm (i)}] ${\mathbb C}P^{n}$ equipped with a symmetric
Fubini-Study metric, or \item[{\rm (ii)}] ${\mathbb C}P^{n}$
equipped with a standard $\Lie{Sp}(m)$-homogeneous metric, $n = 2m-1,$ $m\geq 2,$
or \item[{\rm (iii)}] $S^{6}$ with a Riemannian metric of constant
sectional curvature.
\end{enumerate}
\end{theorem}

\noindent In the proof we make use
of the theory of $3$-symmetric spaces developed by
Gray and Wolf \cite{G,WG} and the
more recent study on homogeneous nearly-K\"ahler manifolds due to
Nagy \cite{Nagy1,Nagy2} and
Butruille \cite{Butrui}.


\section{Preliminaries}\indent

\setcounter{equation}{0}

Let $(M,g)$ be a connected homogeneous {\em Riemannian} manifold.
Then $(M,g)$ can be expressed as
a coset space $G/K,$ where $G$ is a Lie group, which
is assumed to be connected, acting transitively and effectively on $M,$ $K$
is the isotropy subgroup of $G$ at some point $o\in M,$ the origin
of $G/K,$ and $g$ is a $G$-invariant Riemannian metric. Moreover,
we can assume that $G/K$ is a {\em reductive homogeneous space},
i.e. there is an $\mathrm{Ad}(K)$-invariant subspace $\mathfrak{m}$ of the Lie algebra $\mathfrak{g}$ of $G$
such that $\mathfrak{g} = \mathfrak{m} \oplus \mathfrak{k},$
being $\mathfrak{k}$ the Lie algebra of $K.$ $(M=G/K,g)$ is said to be {\em naturally
reductive}, or more precisely $G$-{\em naturally reductive}, if
there exists a reductive decomposition $\mathfrak{g} =
\mathfrak{m} \oplus \mathfrak{k}$ satisfying
\begin{equation}\label{nred}
\langle [X,Y]_\mathfrak{m},Z
\rangle + \langle [X,Z]_\mathfrak{m},Y \rangle = 0,
\end{equation}
\noindent for all $X,Y,Z\in \mathfrak{m},$ where
$[X,Y]_\mathfrak{m}$ denotes the $\mathfrak{m}$-component of $[X,Y]$
and $\langle \cdot , \cdot \rangle$
is the inner product induced by $g$ on
$\mathfrak{m}$,  using the
canonical identification $\mathfrak{m}\cong
\mathrm{T}_{o}M.$ When there exists
a bi-invariant inner product $q$ on ${\mathfrak g}$ whose restriction to
${\mathfrak m} = {\mathfrak k}^{\bot}$ is $\langle \cdot , \cdot \rangle$,
the homogeneous Riemannian manifold $(M=G/K,g)$ is called {\em
normal homogeneous}. Then, for all $X,Y,Z\in {\mathfrak g},$ we have
\begin{equation}\label{tt2}
q([X,Y],Z) + q([X,Z],Y) = 0.
\end{equation}
Hence each normal homogeneous
space is naturally reductive. It is well known that there exists a
bi-invariant inner product on the Lie algebra ${\mathfrak g}$ of
the Lie group $G$ if and only if $G$ is compact. Thus every normal
homogeneous space is compact. Since $G$ is compact and semisimple
if and only if the Cartan-Killing form $B$ is negative definite,
we may choose $q = -B,$ in which case $(M =G/K,g)$ is called a
{\em standard Riemannian homogeneous space} and the induced
Riemannian metric $g$ is called the {\em standard homogeneous
metric} or {\em Cartan-Killing metric} on $M.$

If $G$ is a simple compact Lie group, any naturally reductive $G$-homogeneous
Riemannian manifold is standard, up to scaling factor. Moreover, the unique $G$-invariant
 Riemannian metric, up to homotheties, on a compact isotropy irreducible space $M = G/K$
 is standard choosing the appropriate scaling factor and this Riemannian metric is Einstein.
 Note that not all standard homogeneous metrics are Einstein.

The sectional curvature of a normal homogeneous Riemannian manifold $(M = G/K,g)$ is given by
\begin{equation}\label{curnor}
\langle R(X,Y)X,Y \rangle =
\|[X,Y]_\mathfrak{k}\|^{2} + \frac{\textstyle 1}{\textstyle
4}\|[X,Y]_\mathfrak{m}\|^{2},
\end{equation}
for all $X,Y\in \mathfrak{m}\cong
\mathrm{T}_{o}M.$ So the sectional curvature of a
normal homogeneous manifold is always non-negative and there exists
a section $\pi=\mathbb{R}\{X,Y\},$ $X,Y\in {\mathfrak m},$ such that
$K(\pi) = 0$ if and only if $[X,Y] = 0.$

We shall also need some general results of complex simple Lie
algebras. See \cite{He} for more details. Let ${\mathfrak
g}_{\mathbb C}$ be a simple Lie algebra over ${\mathbb C}$ and
${\mathfrak h}_{\mathbb C}$ a Cartan subalgebra of ${\mathfrak
g}_{\mathbb C}.$ Let $\Delta$ denote the set of non-zero roots of
${\mathfrak g}_{\mathbb C}$ with respect to ${\mathfrak h}_{\mathbb
C}$ and $\Pi = \{\alpha_{1},\dots ,\alpha_{l}\}$ a system of simple
roots or a basis of $\Delta.$ Because the restriction of the
Cartan-Killing form $B$ of ${\mathfrak g}_{\mathbb C}$ to
${\mathfrak h}_{\mathbb C} \times {\mathfrak h}_{\mathbb C}$ is
non-degenerate, there exists a
unique element $H_{\alpha}\in {\mathfrak h}_{\mathbb C}$ such that
$B(H,H_{\alpha}) = \alpha(H),$ for all $H\in {\mathfrak h}_{\mathbb
C}.$ Moreover, we have ${\mathfrak h}_{\mathbb C} = \sum_{\alpha\in
\Delta}{\mathbb C}H_{\alpha}$ and $B$ is strictly positive definite
on ${\mathfrak h}_{\mathbb{R}} = \sum_{\alpha\in \Delta}\mathbb{R}
H_{\alpha}.$ Put $\langle \alpha
,\beta \rangle = B(H_{\alpha},H_{\beta}).$ We choose root vectors
$\{E_{\alpha}\}_{\alpha \in \Delta},$ such that for all
$\alpha,\beta\in \Delta,$ we have
\begin{equation}\label{v}
\left.
\begin{array}{lcl}
[E_{\alpha},E_{-\alpha}]= H_{\alpha},& & [H,E_{\alpha}] =
\alpha(H)E_{\alpha},\;\;\;\mbox{for}\;H\in
{\mathfrak h}_{\mathbb
C};\\[0.6pc]
[E_{\alpha},E_{\beta}] = 0 , &
&\mbox{if}\; \alpha+\beta\neq
0\;\mbox{and}\;\alpha + \beta \not\in \Delta;\\[0.6pc]
[E_{\alpha},E_{\beta}] = N_{\alpha,\beta}E_{\alpha + \beta}
, & &\mbox{if}\;\alpha + \beta
\in \Delta,
\end{array}
\right\}
\end{equation}
where the constants $N_{\alpha,\beta}$ satisfy $N_{\alpha,\beta} = -N_{-\alpha,-\beta},$ $N_{\alpha,\beta} =
-N_{\beta,\alpha}.$ Moreover, given an $\alpha$-series $\beta + n\alpha$ $(p\leq n\leq
q)$ containing $\beta,$ then
\begin{equation}\label{***}
(N_{\alpha,\beta})^{2} = \frac{\textstyle q(1-p)}{\textstyle 2}
  \langle \alpha,\alpha \rangle.
\end{equation}
For this choice, if $\alpha + \beta
\neq 0$, then  $ E_{\alpha}$ and $E_{\beta}$ are orthogonal under
$B$, $B(E_{\alpha},E_{-\alpha}) = 1$ and we have the orthogonal
direct sum
\[
{\mathfrak g}_{\mathbb C} = {\mathfrak h}_{\mathbb C} +
\displaystyle\sum_{\alpha\in \Delta} {\mathbb C}E_{\alpha}.
\]
Denote by $\Delta^{+}$ the set of positive roots of $\Delta$ with
respect to some lexicographic order in $\Pi.$ Then the
$\mathbb{R}$-linear subspace ${\mathfrak g}$ of ${\mathfrak
g}_{\mathbb C}$ given by
\[
{\mathfrak g} = {\mathfrak h} + \displaystyle\sum_{\alpha \in
\Delta^{+}} (\mathbb{R}\; U^{0}_{\alpha} + \mathbb{R}\;
U^{1}_{\alpha})
\]
is a compact real form of ${\mathfrak g}_{\mathbb C},$ where
${\mathfrak h} = \sum_{\alpha\in \Delta}\mathbb{R}
\sqrt{-1}H_{\alpha}$ and $U^{0}_{\alpha} = E_{\alpha}-E_{-\alpha}$
and $U^{1}_{\alpha} = \sqrt{-1}(E_{\alpha} + E_{-\alpha}).$ Next
we put $N_{\alpha,\beta} = 0$ if $\alpha + \beta \not\in \Delta.$
Then, using (\ref{v}), one gets
\begin{lemma}\label{bracket} For all $\alpha ,\beta \in
\Delta^{+}$ and $a = 0,1,$ the following equalities hold:
\begin{enumerate}
\item[{\rm (i)}]
 $[U^{a}_{\alpha},\sqrt{-1} H_{\beta}] = (-1)^{a+1}
 \langle \alpha,\beta \rangle
  U^{a+1}_{\alpha};$ \item[{\rm (ii)}]
$[U^{0}_{\alpha},U^{1}_{\alpha}] = 2\sqrt{-1}H_{\alpha};$
\item[{\rm (iii)}] $[U^{a}_{\alpha},U^{b}_{\beta}] =
(-1)^{ab}N_{\alpha,\beta}U^{a+b}_{\alpha + \beta} +
(-1)^{a+b}N_{-\alpha,\beta}U^{a+b}_{\alpha -\beta},$ where $\alpha
\neq \beta$ and $a\leq b.$
\end{enumerate}
\end{lemma}

For each $\sqrt{-1}H\in {\mathfrak h},$ it implies that
\begin{equation}\label{adjoint}
\begin{array}{lcl}
\mathrm{Ad}_{\exp\sqrt{-1}H}U^{0}_{\alpha}
& = & \cos\alpha(H) U^{0}_{\alpha} +
\sin\alpha(H) U^{1}_{\alpha},\\[0.5pc]
\mathrm{Ad}_{\exp\sqrt{-1}H}U^{1}_{\alpha}
& = & \cos\alpha(H) U^{1}_{\alpha} - \sin\alpha(H) U^{0}_{\alpha}.
\end{array}
\end{equation}

\section{Compact 3-symmetric spaces with strictly positive curvature}\indent
\label{3symmetrcispaces}

\setcounter{equation}{0}

We recall that a connected Riemannian manifold $(M,g)$ is called a
$3$-{\em symmetric space} \cite{G} if it admits a family of
isometries $\{\theta_{p}\}_{p\in M}$ of $(M,g)$ satisfying
\begin{enumerate}
\item[{\rm (i)}] $\theta^{3}_{p} = I,$
\item[{\rm (ii)}] $p$ is an isolated fixed point of $\theta_{p},$
\item[{\rm (iii)}] the tensor field $\Theta$ defined by $\Theta =
(\theta_{p})_{*p}$ is of class $C^{\infty},$
\item[{\rm (iv)}] $\theta_{p*}\comp J = J \comp \theta_{p*},$
\end{enumerate}
where $J$ is the {\em canonical almost complex structure} associated
with the family $\{\theta_{p}\}_{p\in M}$ given by $J =
\frac{1}{\sqrt{3}}(2\Theta +I).$ Riemannian $3$-symmetric spaces are
characterised by a triple
$(G/K,\sigma, \langle \cdot ,
\cdot\rangle)$ satisfying the following conditions:
\begin{enumerate}
\item[{\rm (1)}] $G$ is a connected Lie group and $\sigma$ is an
automorphism of $G$ of order $3,$
 \item[{\rm (2)}] $K$ is a closed
subgroup of $G$ such that $G_{o}^{\sigma} \subseteq K\subseteq
G^{\sigma},$ where $G^{\sigma} = \{x\in G\mid \sigma(x) = x\}$ and
$G_{o}^{\sigma}$ denotes its identity component,
 \item[{\rm (3)}]
$\langle \cdot , \cdot\rangle$ is
an $\mathrm{Ad}(K)$- and
$\sigma$-invariant inner product on the vector space ${\mathfrak
m} = ({\mathfrak m}^{+}\oplus {\mathfrak m}^{-})\cap {\mathfrak
g},$ where ${\mathfrak m}^{+}$ and ${\mathfrak m}^{-}$ are the
eigenspaces of $\sigma$ on the complexification ${\mathfrak
g}_{\mathbb C}$ of ${\mathfrak g}$ corresponding to the
eigenvalues $\varepsilon$ and $\varepsilon^{2},$ respectively,
where $\varepsilon = e^{2\pi\sqrt{-1}/3}.$
\end{enumerate}
Here and in the sequel, $\sigma$ and its differential $\sigma_{*}$
on ${\mathfrak g}$ and on ${\mathfrak g}_{\mathbb C}$ are denoted
by the same letter $\sigma.$ The inner product
$\langle \cdot , \cdot\rangle$
induces a $G$-invariant Riemannian metric $g$ on $M = G/K$ and
$(G/K,g)$ becomes into a Riemannian $3$-symmetric space.
Then it is a reductive
homogeneous space with reductive decomposition ${\mathfrak g} =
{\mathfrak m}\oplus {\mathfrak k},$ where the algebra of Lie
${\mathfrak k}$ of $K$ is ${\mathfrak g}^{\sigma} = \{ X\in
{\mathfrak g}\mid \sigma X = X\}.$ The canonical almost structure
$J$ on $G/K$ is $G$-invariant and it is determined by the
$\mathrm{Ad}(K)$-invariant
automorphism $J_{o}$ on ${\mathfrak m}$ given by
\begin{equation}\label{J}
J_{o} = \frac{1}{\sqrt{3}}(2\sigma_{\mid {\mathfrak m}} +
Id_{\mathfrak m}).
\end{equation}
Note that, taking into account that $\sigma^{2}X + \sigma X + X \in
{\mathfrak k},$ one obtains that in fact $J_{o}^{2} = -Id_{\mathfrak
m}.$ Moreover $(M = G/K,g,J)$ is quasi-K\"ahlerian and it is nearly
K\"ahlerian if and only if $(G/K,g)$ is a naturally reductive
homogeneous space with adapted reductive decomposition ${\mathfrak
g} = {\mathfrak m} \oplus {\mathfrak k}.$ In this case $g$ is said
to be an {\em adapted naturally reductive metric} for $M.$ Under the
canonical identification of ${\mathfrak m}$ with
$\mathrm{T}_{o}G/K,$ we have the
following (see \cite{G}).
\begin{equation}\label{JJ}
[JX,JY]_{\mathfrak k} = [X,Y]_{\mathfrak k},\;\;\;\;\;
[JX,Y]_{\mathfrak m} = -J[X,Y]_{\mathfrak m}.
\end{equation}

According to Gray \cite{G}, a simply connected Riemannian
$3$-symmetric space $(M,g)$
may be decomposed as a Riemannian product
$M = M_{0}\times M_{1}\times \dots \times M_{r},$ where $M_{0}$ is
an even dimensional Euclidean space and $M_{1},\dots , M_{r}$ are
irreducible Riemannian $3$-symmetric spaces. Each $M_{i},$ $i =1,\dots ,r,$ admits a homogeneous metric $g,$ unique up to a scalar multiple, that is nearly K\"ahler and makes $(M_{i},g)$ a standard naturally reductive homogeneous space. A compact irreducible
Riemannian $3$-symmetric space $(M =
G/K,\sigma, \langle \cdot ,
\cdot\rangle)$ has one of the following forms:

\noindent ${\sf Type \;A_{3}:}$ $G$ is a compact connected simple
Lie group acting effectively and $\sigma$ is an inner automorphism
on the Lie algebra ${\mathfrak g}$ of $G.$

Let $\mu = \sum_{i=1}^{l}m_{i}\alpha_{i}$ be the {\em maximal
root} of $\Delta$ and consider $H_{i}\in {\mathfrak h}_{\mathbb
C},$ $i=1,\dots ,l,$ defined by
\[ \alpha_{j}(H_{i}) = \frac{\textstyle 1}{\textstyle m_{i}}\delta_{ij},\;\;\; i,j = 1,\dots ,l.
\]
Following \cite[Theorem 3.3]{WG}, each inner automorphism of order
$3$ on ${\mathfrak g}_{\mathbb C}$ is conjugate in the inner
automorphism group of ${\mathfrak g}_{\mathbb C}$ to some
$\sigma = \mathrm{Ad}_{\exp 2\pi\sqrt{-1}H},$
 where $H = \frac{1}{3}m_{i}H_{i}$ with $1\leq
m_{i}\leq 3$ or $H = \frac{1}{3}(H_{i} + H_{j})$ with $m_{i} =
m_{j} = 1.$ Then there are four classes of $\sigma =
\mathrm{Ad}_{\exp 2\pi\sqrt{-1}H}$ with corresponding simple root
systems $\Pi(H)$ for ${\mathfrak g}^{\sigma}_{\mathbb C},$ Types
$A_{3}I$-$A_{3}IV$ given in Table I. Denote by $\Delta^{+}(H)$ the
positive root system generated by $\Pi(H)$.
Then we have ${\mathfrak
h}\subset {\mathfrak k} = {\mathfrak g}^{\sigma}$ and
\begin{table}[tbp]
{\small
$$
\begin{array}{|c|l|c|c|}
\hline \mbox{\rm Type} & \qquad \quad  \sigma & m_{i} & \Pi(H)\\[0.6pc]\hline
A_{3}I &
\mathrm{Ad}_{\exp\frac{2\pi\sqrt{-1}}{3}H_{i}}& 1
& \{\alpha_{k}\in \Pi\mid k\neq i\}\\[0.3pc]\hline
A_{3} II &
\mathrm{Ad}_{\exp 2\pi\sqrt{-1}\frac{(H_{i}+ H_{j})}{3}} & m_{i} = m_{j} =1
 & \{\alpha_{k}\in \Pi\mid k\neq i,\;k\neq j\}\\[0.3pc]\hline
A_{3} III &
\mathrm{Ad}_{\exp\frac{4\pi\sqrt{-1}}{3}H_{i}}
&
2 & \{\alpha_{k}\in \Pi\mid k\neq i\}\\[0.3pc]\hline
A_{3} IV & \mathrm{Ad}_{\exp
2\pi\sqrt{-1}H_{i}} &
3 & \{\alpha_{k}\in \Pi\mid k\neq i\} \cup \{-\mu\}\\[0.3pc]\hline
\end{array}
$$
} \center{\small Table I}
\end{table}
\[
{\mathfrak k} = {\mathfrak h} + \displaystyle\sum_{\alpha \in
\Delta^{+}(H)}(\mathbb{R}\; U^{0}_{\alpha}+ \mathbb{R}\;
U^{1}_{\alpha}).
\]
Because $B(U^{a}_{\alpha},U^{b}_{\beta}) =
-2\delta_{\alpha\beta}\delta_{ab},$ it follows that
$\{U^{a}_{\alpha}\mid a=0,1,\; \alpha\in
\Delta^{+}\smallsetminus\Delta^{+}(H)\}$ becomes into an
orthonormal basis for $({\mathfrak
m},\langle \cdot , \cdot \rangle=
-\frac{\textstyle 1}{\textstyle 2}B_{\mid {\mathfrak m}}).$
\vspace{0.1cm}

\noindent ${\sf Type\; B_{3}:}$ $G$ is a compact simple Lie group
and the complexification $\mathfrak{g}_\mathbb{C}$ of $\mathfrak{g}$
is of Dynkin type
$\mathfrak{d}_{4}$ and $\sigma$ is an outer automorphism on
$\mathfrak{g}.$

\noindent ${\sf Type\; C_{3}:}$ $G = L\times L\times L,$ where $L$
is a compact simple Lie group and $\sigma$ on ${\mathfrak g} =
{\mathfrak l}\oplus {\mathfrak l}\oplus {\mathfrak l}$ is given by
$\sigma(X,Y,Z) = (Z,X,Y).$ Here, ${\mathfrak k} = {\mathfrak
g}^{\sigma}$ is ${\mathfrak l}$ embedded diagonally.

\begin{lemma}\label{linner} If $\sigma$ is an inner automorphism,
i.e. it is of Type $A_{3},$ then $J_{o}$ defined as in
{\rm (\ref{J})} satisfies
\[
J_{o}U^{0}_{\alpha} = \pm U^{1}_{\alpha} ,\;\;\;\; JU^{1}_{\alpha}
= \mp U^{0}_{\alpha},\;\;\;\mbox{for
all}\;\;\alpha\in\Delta^{+}\smallsetminus\Delta^{+}(H).
\]
\end{lemma}
{\sf Proof.} Each $\alpha\in\Delta^{+}\smallsetminus\Delta^{+}(H)$
may be written as $ \alpha = \sum_{k=1}^{l}n_{k}\alpha_{k}$, where
$n_{k}\in{\mathbb Z},$ $n_{k}\geq 0,$ for each $k\in\{1,\dots,
l\}.$ Let $\sqrt{-1}H\in {\mathfrak h}$ such that $\sigma =
\mathrm{Ad}_{\exp 2\pi\sqrt{-1}H},$ as before. Then we have
$\alpha(H) = \frac13 n_{i},$ for $H=\frac{1}{3}m_{i}H_{i},$ $1\leq
m_{i}\leq 3,$ and $\alpha(H)  = \frac{1}{3}(n_{i} + n_{j}),$ for
$H =\frac{1}{3}(H_{i} + H_{j}),$ $m_{i} = m_{j} =1.$ Hence it
follows that the possible values of $\alpha(H)$ are $\frac{1}{3},$
$\frac{2}{3}$ and $1.$ But $\alpha(H)\neq 1,$ because for
$\alpha(H) =1,$ one obtains from (\ref{adjoint}) that
$\sigma(U^{0}_{\alpha}) = U^{0}_{\alpha}.$ Then, the result
follows directly using again (\ref{adjoint}). \hfill $\Box$

In \cite{Wa} Wallach showed  that
a simply connected $2n$-dimensional, compact homogeneous
Riemannian manifold with strictly positive sectional curvature is
isometric to
\begin{enumerate}
\item[{\rm (i)}] a compact rank one symmetric space: ${\mathbb
C}P^{n},$ $S^{2n},$ ${\mathbb H}P^{n/2}$ $(n$ even{\rm )},
${\mathbb C}aP^{2}$ $(n =8);$

\noindent or

\item[{\rm (ii)}] one of the following quotient spaces $G/K$ with
a suitable $G$-invariant metric:
\begin{enumerate}
\item[{\rm (1)}] the manifolds of flags in the complex,
quaternionic and Cayley three-space:
$$
\begin{array}{lcl}
{\mathbb F}^{6} & = & \Lie{SU}(3)/(\Lie{U}(1)\times \Lie{U}(1)),\\
{\mathbb F}^{12} & = & \Lie{Sp}(3)/(\Lie{SU}(2)\times \Lie{SU}(2)\times \Lie{SU}(2)),\\
{\mathbb F}^{24} & = & F_{4}/\Lie{Spin}(8);
\end{array}
$$
\item[{\rm (2)}] $\Lie{Sp}(m)/(\Lie{Sp}(m-1)\times \Lie{U}(1)),$
$n = 2m-1;$ \item[{\rm (3)}] the six-dimensional sphere $S^{6} =
\Lie{G}_{2}/\Lie{SU}(3).$
\end{enumerate}
\end{enumerate}

\noindent Then we have
\begin{proposition}\label{p3-sym} A simply connected, compact
Riemannian $3$-symmetric space with strictly positive sectional
curvature is isometric to one of the following spaces, equipped
with a suitable invariant metric:
\begin{enumerate}
 \item[{\rm (i)}] $A_{3}I:$ ${\mathbb C}P^{n} =
\Lie{SU}(n+1)/\Lie{S}(\Lie{U}(1)\times \Lie{U}(n));$
  \item[{\rm (ii)}] $A_{3}II:$ ${\mathbb F}^{6} = \Lie{SU}(3)/(\Lie{U}(1)\times
\Lie{U}(1));$
 \item[{\rm (iii)}] $A_{3}III:$ ${\mathbb C}P^{n} =
\Lie{Sp}(m)/(\Lie{Sp}(m-1)\times \Lie{U}(1)),$ $n = 2m-1;$
 \item[{\rm (iv)}] $A_{3}IV:$ $S^{6} = \Lie{G}_{2}/\Lie{SU}(3).$
\end{enumerate}

\end{proposition}
\noindent{\sf Proof.} A compact irreducible Riemannian $3$-symmetric
space $(M = G/K,\sigma,\langle \cdot
, \cdot \rangle )$ is a (Hermitian) symmetric coset space if and
only if $\sigma$ is of Type $A_{3}I$ (see for example \cite{GD}).
Then if it is also of rank one it
must be the complex projective space ${\mathbb C}P^{n} =
\Lie{SU}(n+1)/\Lie{S}(\Lie{U}(1)\times \Lie{U}(n)).$

Next we show that that the quotient
spaces ${\mathbb F}^{6}= \Lie{SU}(3)/T^{2},$ ${\mathbb C}P^{n}
=\Lie{Sp}(m)/(\Lie{Sp}(m-1)\times \Lie{U}(1))$ and $S^{6} =
\Lie{G}_{2}/\Lie{SU}(3)$ admit a structure of Riemannian
$3$-symmetric space. In $\mathfrak{a}_{2}:\xymatrix@R=.5cm@C=.8cm{
\stackrel{1}{\stackrel{\circ}{\alpha_{1}}} \ar@{-}[r] &
\stackrel{1}{\stackrel{\circ}{\alpha_{2}}}}$ we consider the inner
automorphism $\sigma = \mathrm{Ad}_{\exp 2\pi\sqrt{-1}H}$ of Type
$A_{3}II$ with $H=\frac{1}{3}(H_{1} + H_{2}).$ From Table I, the
simple root system $\Pi(H)$ for ${\mathfrak a}_{2}^{\sigma}$ is
empty and ${\mathfrak a}_{2}^{\sigma}$ is the $2$-dimensional torus
generates by $\{\sqrt{-1}H_{\alpha_{1}}, \sqrt{-1}H_{\alpha_{2}}\}.$

For the complex Lie algebra ${\mathfrak
g}_{\mathbb C}  = {\mathfrak c}_{m}:\xymatrix@R=.5cm@C=.8cm{
\stackrel{2}{\stackrel{\circ}{\alpha_{1}}} \ar@{-}[r] &
\stackrel{2}{\stackrel{\circ}{\alpha_{2}}} \ar@{-}[r] & \; \dots
\ar@{-}[r] & \stackrel{2}{\stackrel{\circ}{\alpha_{m-1}}}
\ar@2{-}[r] & \stackrel{1}{\stackrel{\circ}{\alpha_{m}}}},$ $m\geq
2,$ a set of positive
roots is given by
$$
\begin{array}{lcl}
\Delta^{+} & = & \{\alpha_{ij} = \alpha_{i} + \dots
+\alpha_{j}\;\;
(1\leq i\leq j\leq m),\\[0.5pc]
& & \hspace{0.3cm} \widetilde{\alpha_{ij}} = \alpha_{i} +\dots +
2\alpha_{j} + \dots + 2\alpha_{m-1} + \alpha_{m}\;\; (1\leq i\leq
j\leq m-1)\}.
\end{array}
$$
The inner automorphism $\sigma = \mathrm{Ad}_{\exp
2\pi\sqrt{-1}H},$ where $H = \frac{2}{3}H_{1},$  is of Type
$A_{3}III.$  Then, using Table I, the simple root system $\Pi(H)$
for the (complex) Lie algebra ${\mathfrak c}_{m}^{\sigma} = \{X\in
{\mathfrak c}_{m}\mid \sigma X = X\}$ is $\pi(H) = \{\alpha_{2},
\dots, \alpha_{m}\}.$ The positive root system $\Delta^{+}(H)$
generated by $\Pi(H)$ is
\[
\Delta^{+}(H) =  \{\alpha_{ij}\; (2\leq i\leq j\leq m),\;
\widetilde{\alpha_{ij}}\; (2\leq i\leq j\leq m-1)\}.
\]
Thus $\{H_{\alpha_{i}}\; (1\leq i\leq m); \;E_{\pm \alpha_{ij}} \;
(2\leq i\leq j\leq m),\; E_{\pm (\widetilde{\alpha_{ij}})}\;
(2\leq i\leq j\leq m-1)\}$ is a Weyl basis for ${\mathfrak
c}_{m}^{\sigma}.$ It implies that ${\mathfrak c}_{m}^{\sigma}$ is
of type ${\mathfrak c}_{m-1}\otimes {\mathfrak T}^{1}$ and
${\mathfrak k} = {\mathfrak s}{\mathfrak p}(m-1)\oplus {\mathfrak
u}(1).$

For the exceptional Lie algebra ${\mathfrak
g}_{2}:\;\xymatrix@R=.5cm@C=.3cm{
\stackrel{3}{\stackrel{\circ}{\alpha_{1}}} \ar@3{-}[r] &
\stackrel{2}{\stackrel{\circ}{\alpha_{2}}}},$ a set of positive
roots is given by $\Delta^{+} = \{\alpha_{1},\alpha_{2},\alpha_{1}
+ \alpha_{2}, 2\alpha_{1} + \alpha_{2}, 3\alpha_{1} + \alpha_{2},
3\alpha_{1} + 2\alpha_{2}\},$ being $\mu = 3\alpha_{1} +
2\alpha_{2}$ the maximal root. Next let $\sigma = \mathrm{Ad}_{\exp 2\pi\sqrt{-1}H}$ be the
inner automorphism of ${\mathfrak g}_{2}$ such that $H= H_{1}.$
Then $\sigma$ is of Type
$A_{3}IV$ and, from Table I, the
simple root system $\Pi(H)$ for the (complex) Lie algebra
${\mathfrak g}_{2}^{\sigma} = \{X\in {\mathfrak g}_{2}\mid \sigma
X = X\}$ is $\pi(H) = \{\alpha_{2},-\mu\}.$  Hence the positive
root system $\Delta^{+}(H)$ generated by $\Pi(H)$ is given by
$\Delta^{+}(H) = \{\alpha_{2},3\alpha_{1} + \alpha_{2},\mu\}$ and
$\{H_{\alpha_{1}}, H_{\alpha_{2}}; E_{\pm \alpha_{1}}, E_{\pm
(3\alpha_{1} + \alpha_{2})}, E_{\pm \mu}\}$ is a Weyl basis for
${\mathfrak g}_{2}^{\sigma}.$ It implies that ${\mathfrak
g}_{2}^{\sigma}$ is of type ${\mathfrak a}_{2}$ and
its corresponding compact real
form ${\mathfrak k}$ is isomorphic to ${\mathfrak s}{\mathfrak
u}(3).$

Finally, we prove that the quotient spaces ${\mathbb F}^{12} = \Lie{Sp}(3)/(\Lie{SU}(2)\times \Lie{SU}(2)\times \Lie{SU}(2))$ and ${\mathbb F}^{24} = F_{4}/\Lie{Spin}(8)$ do not
admit any structure of Riemannian $3$-symmetric space. Because any automorphism of order $3$ on ${\mathfrak c}_{3}$ or on ${\mathfrak f}_{4}$ must be of Type $A_{3},$ it can be written , up to conjugance, as $\sigma = \mathrm{Ad}_{\exp 2\pi\sqrt{-1}H}.$ On ${\mathfrak
c}_{3}:\xymatrix@R=.5cm@C=.8cm{
\stackrel{2}{\stackrel{\circ}{\alpha_{1}}} \ar@{-}[r] &
\stackrel{2}{\stackrel{\circ}{\alpha_{2}}} \ar@2{-}[r] &
\stackrel{1}{\stackrel{\circ}{\alpha_{3}}}},$ $H$
is given by $\frac{1}{3}H_{3},$ $\frac{2}{3}H_{1}$ or $\frac{2}{3}H_{2}.$
Then ${\mathfrak c}^{\sigma}_{3}$ is of type ${\mathfrak
a}_{2}\oplus {\mathfrak T}^{1},$ ${\mathfrak c}_{2}\oplus
{\mathfrak T}^{1}$ or $({\mathfrak a}_{1} \oplus {\mathfrak T}^{1})
\oplus {\mathfrak a}_{1},$ respectively. On ${\mathfrak f}_{4}:\xymatrix@R=.5cm@C=.8cm{
\stackrel{2}{\stackrel{\circ}{\alpha_{1}}} \ar@{-}[r] &
\stackrel{3}{\stackrel{\circ}{\alpha_{2}}}\ar@2{-}[r] &
\stackrel{4}{\stackrel{\circ}{\alpha_{3}}}\ar@{-}[r] &
\stackrel{2}{\stackrel{\circ}{\alpha_{4}}}}$
there are three possibilities for $H:$ $\frac{2}{3}H_{1},$
$H_{3}$ or $\frac{2}{3}H_{4}$ and their corresponding complex Lie
algebras ${\mathfrak f}^{\sigma}_{4}$ are of type ${\mathfrak
b}_{3}\oplus {\mathfrak T}^{1},$ ${\mathfrak a}_{2}\oplus
{\mathfrak a}_{2}$ or ${\mathfrak c}_{3}\oplus {\mathfrak T}^{1}$. It gives the desired result.
\hfill
$\Box$


\section{Proof of Theorems 1.1 and 1.2}\indent
\setcounter{equation}{0}

In \cite{Va} Valiev  determined
the set of all homogeneous Riemannian metrics on ${\mathbb
F}^{6},$ ${\mathbb F}^{12}$ and ${\mathbb F}^{24}$ of strictly
positive sectional curvature and their corresponding optimal
pinching constants. According to Berger's classification, they
cannot be normal homogeneous. It is worthwhile to note that $S^{6}
= \Lie{G}_{2}/\Lie{SU}(3)$ carries the usual metric of constant
sectional curvature because $\Lie{SU}(3)$ is irreducible on the
tangent space (see \cite{Be2}). Nevertheless,
$(\Lie{G}_{2},\Lie{SU}(3))$ is not a symmetric pair.

Next, in order to prove Theorem 1.1, we only need to study the
complex projective space ${\mathbb C}P^{n} =
\Lie{Sp}(m)/(\Lie{Sp}(m-1)\times \Lie{U}(1)),$ $n = 2m-1,$
as a standard
$\Lie{Sp}(m)$-homogeneous
 Riemannian manifold. See \cite{Z} for a general study about all its homogeneous Riemannian metrics.
  On $\mathfrak{c}_{m},$ the maximal root $\mu$ is given by $\mu= \widetilde{{\alpha}_{11}}$ and we have
$$
\langle \alpha_{i},\alpha_{i}\rangle = \left\{
\begin{array}{ll}
\tfrac{1}{2(m+1)}, & 1\leq i \leq m-1,\\[0.4pc]
\tfrac{1}{m+1}, & i=m,
\end{array}
\right. \qquad
 \langle \alpha_{j},\alpha_{j+1}\rangle = \left\{
\begin{array}{ll}
-\tfrac{1}{4(m+1)}, & 1\leq j\leq m-2,\\[0.4pc]
-\tfrac{1}{2(m+1)}, & j = m-1
\end{array}
\right.
$$
and
 $\langle \alpha_{i},\alpha_{j}\rangle =
 0$
 for the remaining $(i,j).$ From
Proposition \ref{p3-sym},
$Sp(m)/(Sp(m-1)\times U(1))$ is
an irreducible $3$-symmetric space of Type $A_{3}III.$ A basis for
${\mathfrak k}= {\mathfrak s}{\mathfrak p}(m-1)\oplus {\mathfrak
T}^{1}$ is given by
\[
\{\sqrt{-1}H_{\alpha_{i}}\; (1\leq i\leq m); \;
U^{a}_{\alpha_{ij}} \; (2\leq i\leq j\leq m),\;
U^{a}_{\widetilde{\alpha_{ij}}} \; (2\leq i\leq j\leq m-1), \; a =
0,1\}.
\]
\noindent Denote by $\beta_{i} = \alpha_{1i}$ $(1\leq i\leq m)$ and
$\beta_{m+j} = \widetilde{\alpha_{1(m-j)}}$ $(1\leq j\leq m-1).$
Then $\{U^{a}_{\beta_{i}},U^{a}_{\beta_{m+j}}\}_{a = 0,1}$ forms
an orthonormal basis for $({\mathfrak m},\langle\cdot,\cdot
\rangle = -\frac{1}{2}B_{\mathfrak m}),$ where ${\mathfrak m}$ is
the orthogonal complement of ${\mathfrak k}$ in ${\mathfrak
s}{\mathfrak p}(m)$ with respect to the Killing form $B.$ From
Lemma \ref{bracket}, one obtains that the subspaces ${\mathfrak
m}_{1} = \mathbb{R}\{U^{a}_{\mu}\}_{a = 0,1}\cong {\mathbb C}$ and
${\mathfrak m}_{2} = {\mathfrak m}^{\bot}\cong {\mathbb H}^{n}$ of
${\mathfrak m}$
 are $\mathrm{Ad}(\Lie{Sp}(m-1)\times \Lie{U}(1))$-invariant.

It is easy to see by a case-by-case check the following.
\begin{lemma} We have:
\begin{enumerate}
\item[{\rm (i)}] For $1\leq i<j\leq m,$ $\beta_{i}-\beta_{j} =
-\alpha_{(i+1)j}$ and $\beta_{i} + \beta_{j}\in \Delta$ if and only
if $(i,j) = (m-1,m).$ Then $\beta_{m-1} + \beta_{m} = \beta_{2m-1} = \mu.$
\item[{\rm (ii)}] For $1\leq i<j\leq m-1,$
$\beta_{m+i} - \beta_{m+j} = -\alpha_{(m-j)(m-i-1)}$ and $\beta_{m+i} +
\beta_{m+j}\not\in \Delta.$ In particular, $\beta_{m+i} -
\beta_{2m-1} = -\beta_{m-i-1},$ for $1\leq i\leq m-2.$ \item[{\rm
(iii)}] For $1\leq i\leq m$ and $1\leq j\leq m-1,$ $\beta_{i} +
\beta_{m+j}\in \Delta$ if and only if $j = m-i-1$ and $1\leq i\leq
m-2.$ Then $\beta_{i} +
\beta_{2m-i-1} = \mu.$ Moreover,
$$
-(\beta_{i} - \beta_{m+j}) = \left\{
\begin{array}{lcl}
\beta_{2m-i-1},&{\mbox if} & 1\leq i\leq m-2,\; j = m-1;\\[0.5pc]
\widetilde{\alpha_{(m-j)(i+1)}}, & {\mbox if}& 2\leq m-j\leq i\leq
m-2;\\[0.5pc]
\alpha_{(m-j)m}, &{\mbox if} & i=m-1,\; 1\leq j\leq m-2;\\[0.5pc]
\beta_{m}, &{\mbox if} & i = j = m-1;\\[0.5pc]
\alpha_{(m-j)(m-1)}, &{\mbox if} & i=m,\; 1\leq j\leq m-2;\\[0.5pc]
\beta_{m}, &{\mbox if} & i=m,\; j = m-1;\\[0.5pc]
\widetilde{\alpha_{(i+1)(m-j)}}, &{\mbox if} & i<m-j.
\end{array}
\right.
$$
\end{enumerate}
\end{lemma}

Hence, using Lemma \ref{bracket} (iii), the brackets
$[U^{a}_{\beta_{k}},U^{a}_{\beta_{l}}],$ $a = 0,1,$ $1\leq k<
l\leq 2m-1,$ are given by
\begin{equation}\label{brc}
\begin{array}{lcl}
[U^{a}_{\beta_{i}},U^{a}_{\beta_{j}}] & = &
-N_{-\beta_{i},\beta_{j}}U^{0}_{\alpha_{(i+1)j}},\;\;\; 1\leq
i<j\leq m,\;\; (i,j)\neq (m-1,m);\\[0.5pc]
[U^{a}_{\beta_{m-1}},U^{a}_{\beta_{m}}] & = &
(-1)^{a}N_{\beta_{m-1},\beta_{m}}U^{0}_{\beta_{2m-1}} -
N_{-\beta_{m-1},\beta_{m}}U^{0}_{\alpha_{m}};\\[0.5pc]
[U^{a}_{\beta_{m+i}},U^{a}_{\beta_{m+j}}] & = &
-N_{-\beta_{m+i},\beta_{m+j}}U^{0}_{\alpha_{(m-j)(m-i-1)}},\;\;\;
1\leq i<j\leq m-2;\\[0.5pc]
[U^{a}_{\beta_{m+i}},U^{a}_{\beta_{2m-1}}] & = &
-N_{-\beta_{m+i},\beta_{2m-1}}U^{0}_{\beta_{m-i-1}},\;\;\; 1\leq
i\leq m-2;\\[0.5pc]
[U^{a}_{\beta_{i}},U^{a}_{\beta_{2m-1}}] & = &
-N_{-\beta_{i},\beta_{2m-1}}U^{0}_{\beta_{2m-i-1}},\;\;\;1\leq
i\leq m-2;\\[0.5pc]
[U^{a}_{\beta_{i}},U^{a}_{\beta_{m+j}}] & = &
-N_{-\beta_{i},\beta_{m-j}}U^{0}_{\widetilde{\alpha_{(m-j)(i+1)}}},\;\;\;2\leq
m-j\leq i\leq m-2;\\[0.5pc]
[U^{a}_{\beta_{m-1}},U^{a}_{\beta_{m+j}}] & = &
-N_{-\beta_{m-1},\beta_{m+j}}U^{0}_{\alpha_{(m-j)m}},\;\;\; 1\leq j\leq m-2;\\[0.5pc]
[U^{a}_{\beta_{m-1}},U^{a}_{\beta_{2m-1}}] & = &
-N_{-\beta_{m-1},\beta_{2m-1}}U^{0}_{\beta_{m}};\\[0.5pc]
[U^{a}_{\beta_{m}},U^{a}_{\beta_{m+j}}] & = &
-N_{-\beta_{m},\beta_{m+j}}U^{0}_{\alpha_{(m-j)(m-1)}},\;\;\;1\leq
j\leq m-2;\\[0.5pc]
[U^{a}_{\beta_{m}},U^{a}_{\beta_{2m-1}}] & = &
-N_{-\beta_{m},\beta_{2m-1}}U^{0}_{\beta_{2m-1}};\\[0.5pc]
[U^{a}_{\beta_{i}},U^{a}_{\beta_{2m-i-1}}] & = &
(-1)^{a}N_{\beta_{i},\beta_{2m-i-1}}U^{0}_{\beta_{2m-1}} -
N_{-\beta_{i},\beta_{m}}U^{0}_{\widetilde{\alpha_{(i+1)(i+1)}}},\;\;\;
1\leq i \leq m-2;\\[0.5pc]
[U^{a}_{\beta_{i}},U^{a}_{\beta_{m+j}}] & = &
-N_{-\beta_{i},\beta_{m+j}}U^{0}_{\widetilde{\alpha_{(i+1)(m-j)}}},\;\;\;
1\leq i<m-j-1.
\end{array}
\end{equation}

\begin{proposition}\label{ccp} The complex projective space ${\mathbb C}P^{n} =
\Lie{Sp}(m)/(\Lie{Sp}(m-1)\times \Lie{U}(1)),$ $n = 2m-1,$ $m\geq 2,$
equipped with the standard $\Lie{Sp}(m)$-homogeneous Riemannian
metric has strictly sectional curvature with pinching constant
$\delta = \frac{1}{16}.$
\end{proposition}
\noindent{\sf Proof.} Put
\[
u = \sum_{k=1 \atop a=0,1}^{2m-1}u_{a}^{k}U^{a}_{\beta_{k}},\;\;\; v
=\sum_{k=1 \atop a=0,1}^{2m-1}v_{a}^{k}U^{a}_{\beta_{k}}
\]
two elements of ${\mathfrak m}$ and suppose that $[u,v] =0.$ From
Lemma \ref{bracket} (ii), we have
\[
\frac{1}{2}[u,v]_{\mathfrak h} =
\sum_{k=1}^{2m-1}M_{(k,0)(k,1)}\sqrt{-1}H_{\beta_{k}},
\]
where $M_{(k,a)(l,b)}= u^{k}_{a}v_{b}^{l}- u^{l}_{b}v^{k}_{a},$
for each  $k,j=1,\dots, 2m-1$ and $a,b=0,1,$ is the minor in the
matrix of coefficients of $u$ and $v$ with respect to the basis
$\{U^{a}_{\beta_{k}}\}.$ The elements $\sqrt{-1}H_{\beta_{i}} =
\sqrt{-1}(H_{\alpha_{1}} + \dots + H_{\alpha_{i}}),$ $1\leq i\leq
m,$ constitute a basis for ${\mathfrak h}$ and, moreover, one gets
\[
H_{\beta_{m+j}} = H_{\beta_{m}} + H_{\beta_{m-1}} -
H_{\beta_{m-j-1}},\;(1\leq j\leq m-2);\;\;\; H_{\beta_{2m-1}} =
H_{\beta_{m}} + H_{\beta_{m-1}}.
\]
Hence it follows
\begin{equation}\label{H1}
\begin{array}{rcl}
M_{(i,0)(i,1)} & = & M_{(2m-i-1,0)(2m-i-1,1)},\;\;\; i = 1,\dots,
m-2;\\[0.5pc]
M_{(m-1,0)(m-1,1)} & = & M_{(m,0)(m,1)} =
-\displaystyle\sum_{j=1}^{m-1}M_{(m+j,0)(m+j,1)}.
\end{array}
\end{equation}
Using again Lemma \ref{bracket}, equation (\ref{brc}) and taking the coefficients of $[u,v]$ in
$U^{0}_{\beta_{m-i-1}},$ $U^{0}_{\beta_{2m-i-1}},$
$U^{0}_{\widetilde{\alpha_{(i+1)(i+1)}}},$ for $1\leq i\leq m-2,$ and
$U^{0}_{\beta_{m-1}},$ $U^{0}_{\beta_{m}}$ and $U^{0}_{\alpha_{m}}$
we have, respectively, the following equations:
$$
\begin{array}{lcl}
M_{(m+i,0)(2m-1,0)} + M_{(m+i,1)(2m-1,1)} & = & 0;\\[0.5pc]
M_{(i,0)(2m-1,0)} + M_{(i,1)(2m-1,1)} & = & 0;\\[0.5pc]
M_{(i,0)(2m-i-1,0)} + M_{(i,1)(2m-i-1,1)} & = & 0,
\end{array}
$$
for $1\leq i\leq m-2,$ and
$$
\begin{array}{lcl}
M_{(m,0)(2m-1,0)} + M_{(m,1)(2m-1,1)} & = & 0;\\[0.5pc]
M_{(m-1,0)(2m-1,0)} + M_{(m-1,1)(2m-1,1)} & = & 0;\\[0.5pc]
M_{(m-1,0)(m,0)} + M_{(m-1,1)(m,1)} & = & 0.
\end{array}
$$
From here and from (\ref{H1}), if $[u,v]=0$ then one obtains
\begin{equation}\label{H2}
\left.
\begin{array}{rcl}
 M_{(i,0)(i,1)} & = & M_{(2m-i-1,0)(2m-i-1,1)},\;\;\; i =
1,\dots,
m-2;\\[0.5pc]
M_{(m-1,0)(m-1,1)} & = & M_{(m,0)(m,1)} =
-\sum_{j=1}^{m-1}M_{(m+j,0)(m+j,1)};\\[0.5pc]
M_{(i,0)(2m-i-1,0)} + M_{(i,1)(2m-i-1,1)} & = & 0,\;\;\; i
=1,\dots m-1;\\[0.5pc]
M_{(k,0)(2m-1,0)} + M_{(k,1)(2m-1,1)} & = & 0,\;\;\; k =1,\dots
2m-2.
\end{array}
\right\}
\end{equation}
This implies, making straight calculations, that $u$ and $v$ are
linearly dependent and so the sectional curvature of the standard
$\Lie{Sp}(m)$-homogeneous Riemannian metric is strictly positive.
Using \eqref{curnor} and Lemma \ref{bracket} (ii), one gets that the
sectional curvature $K({\mathfrak m}_{1})$ is given by
$2\langle \mu,\mu \rangle$. From \eqref{brc}, one also gets
$K(U^{a}_{\beta_{k}},U^{a}_{\mu}) =
\frac{1}{4}N^{2}_{-\beta_{k},\mu},$ $ k =1,\dots, 2(m-1).$
Then, applying \eqref{***}, we can
conclude
\begin{equation}\label{curm1}
K({\mathfrak m}_{1}) = \frac{2}{m+1},\;\;\; K(X,Y) = \frac{1}{8(m+1)},
\;\;\mbox{if}\;X\in {\mathfrak m}_{1},\;\; Y\in {\mathfrak m}_{2}.
\end{equation}
According with \cite[p. 357]{Z} they are the maximum and the
minimum value for  the sectional curvature, respectively.
Hence the pinching of the
$\Lie{Sp}(m)$-standard metric
$\langle \cdot , \cdot \rangle$
is $\delta = \frac{1}{16}.$ \hfill $\Box$
\begin{remark}{\rm The $\Lie{Sp}(m)$-standard metric on ${\mathbb C}P^{n},$ $n = 2m-1,$ $m\geq 2,$ is
not the symmetric Fubini-Study because, as it is well known, the
pinching constant of this last one is $\delta = \frac{1}{4}.$
Moreover, the $\Lie{Sp}(m)$-standard metric corresponds up to a
constant with the metric defined in \cite[p. 356-357]{Z} for $t =
\frac{1}{2},$
 we have that it is Einstein if and only if $n = 3.$
}
\end{remark}
\begin{lemma} Any simply connected, irreducible
non-K\"{a}hler homogeneous nearly K\"{a}hler manifold
$(M,g,J)$ is a compact $3$-symmetric space and $J$ is its
canonical almost complex structure.
\end{lemma}
\noindent{\sf Proof.} From \cite[Proposition 2.1]{Nagy1},
$(M,g,J)$ is strict nearly K\"ahler. Then, using Nagy
\cite{Nagy2} and Butruille \cite{Butrui}, it must be a
$3$-symmetric space. Moreover, from \cite[Theorem 1.1 (ii)]{Nagy1}
it is compact with a finite fundamental group.\hfill $\Box$

From here, we can reduce our study to nearly K\"{a}hler
$3$-symmetric spaces. Since a compact nearly K\"{a}hler manifold
with positive curvature is simply connected \cite{G1}, using
Theorem \ref{mean1} and Proposition \ref{p3-sym}, we have
\begin{proposition}\label{pmean2} Any non-K\"ahler homogeneous nearly K\"{a}hler manifold with strictly
positive sectional curvature is holomorphically isometric to one
of the following $3$-symmetric spaces, with respect to the
canonical complex structure:
\[
A_{3}III:\;{\mathbb C}P^{n} =
\Lie{Sp}(m)/(\Lie{Sp}(m-1)\times \Lie{U}(1)),\;n = 2m-1;\;\;\;\;\;A_{3}IV:\;S^{6} = \Lie{G}_{2}/\Lie{SU}(3).
\]
\end{proposition}
\noindent Finally, since  compact K\"ahler manifolds with strictly
positive sectional curvature and constant scalar curvature are
isometric to a complex projective space with the Fubini-Study
metric (see \cite{G3}), then Theorem \ref{mean2} follows as a
direct consequence from Proposition \ref{pmean2}.

\end{document}